\documentclass[11pt]{amsart}
\usepackage{amssymb, amsmath}
\usepackage{amsthm}
\usepackage{graphicx}
\newtheorem{theorem}{Theorem}
\newtheorem{lemma}[theorem]{Lemma}

\newtheorem{corollary}[theorem]{Corollary}

\theoremstyle{definition}
\newtheorem*{remark}{Remark}
\title{Smooth (non)rigidity of cusp-decomposable manifolds}
\author{T. Tam Nguyen Phan}
\address{Department of Mathematics\\
5734 S. University Ave.\\
Chicago, IL 60637}
\email{ttamnp@math.uchicago.edu}

\DeclareMathOperator{\Ker}{Ker}

\DeclareMathOperator{\Stab}{Stab}
\DeclareMathOperator{\Isom}{Isom}
\DeclareMathOperator{\Out}{Out}
\DeclareMathOperator{\OO}{O}

\DeclareMathOperator{\Aut}{Aut}

\DeclareMathOperator{\Diff}{Diff}
\DeclareMathOperator{\Conj}{Conj}

\def\R{\mathbb{R}}
\def\Z{\mathbb{Z}}
\def\N{\mathbb{N}}
\def\Q{\mathbb{Q}}
\def\C{\mathbb{C}}
\def\A{\mathbb{O}}

\def\h{\mathfrak{h}}
\begin{document}
\begin{abstract}
We define \emph{cusp-decomposable} manifolds and prove smooth rigidity within this class of manifolds. These manifolds generally do not admit a nonpositively curved metric but can be decomposed into pieces that are diffeomorphic to finite volume, locally symmetric, negatively curved manifolds with cusps. We prove that the group of outer automorphisms of the fundamental group of such a manifold is an extension of an abelian group by a finite group. Elements of the abelian group are induced by diffeomorphisms that are analogous to Dehn twists in surface topology. We also prove that the outer automophism group can be realized by a group of diffeomorphisms of the manifold.
\end{abstract}
\maketitle
\section{Introduction}

Let $M$ be a smooth manifold. We denote by $\Diff(M)$ the group of self-diffeomorphisms of $M$ and by $\Diff_0(M)$ the group of self-diffeomorphisms of $M$ that are homotopic to the identity map. Hence, $\Diff(M)/\Diff_0(M)$ is the group of self-diffeomorphisms of $M$ up to homotopy. Let $\Out(\pi_1(M))$ be the group of outer automorphisms of $\pi_1(M)$. The action of a diffeomorphism on $\pi_1(M)$ induces a natural homomorphism 

\begin{equation}\label{*}
\eta \colon \Diff(M)/\Diff_0(M) \longrightarrow \Out(\pi_1(M)).
\end{equation}
If $M$ is aspherical, then $\eta$ is always an injection. 

It is known for a number of classes of manifolds that $\eta$ is an isomorphism. These include closed surfaces, by the Dehn-Nielsen-Baer Theorem; infra-nil manifolds, by Auslander \cite{Auslander}; and finite-volume, complete, irreducible, locally symmetric, nonpositively curved manifold of dimension greater than $2$, by Mostow Rigidity. If we relax the smoothness condition and require only continuity, i.e. replace diffeomorphsims by homeomorphisms, this also holds for closed, nonpositively curved manifolds of dimension greater than $4$, by Farrell and Jones \cite{FJ2}, and for solvmanifolds, by work of Mostow \cite{Mostow}.
 
The purpose of this paper is to prove that the homomorphism $\eta$ is an isomorphism for a class of manifolds that we call \emph{cusp-decomposable}. These are manifolds that are obtained by taking finite-volume, complete, locally symmetric, negatively curved manifolds with deleted cusps and gluing them together along their cuspidal boundaries via affine diffeomorphisms. (See Section~\ref{sec:definitions} for a precise definition of cusp-decomposable manifolds). For example, the double of a finite-volume, complete, locally symmetric, cusped, negatively curved manifold is cusp-decomposable. Any finite cover of such a manifold is also a cusp-decomposable manifold. 

Cusp-decomposable manifolds do not admit locally homogeneous Riemannian metrics \cite{Tam}. If $M$ is cusp-decomposable and all the pieces in the cusp decomposition of $M$ are hyperbolic manifolds, then $M$ admits a nonpositively curved Riemannian metric \cite{AF}. However, if the pieces in the cusp decomposition of $M$ are complex, quarternionic or octonionic hyperbolic, then $M$ does not generally admit any nonpositively curved metric \cite{Tam}. Hence, the above rigidity results do not apply to this class of manifolds

The first main result in this paper is the following. 

\begin{theorem}\label{Mostow}

Let $M$ and $N$ be cusp-decomposable manifolds of dimension greater than $2$. Then any isomorphism of fundamental groups $\pi_1(M) \cong \pi_1(N)$ is realized by a diffeomorphism from $M$ to $N$. Therefore, there is an isomorphism
\[ \rho \colon \Out(\pi_1(M)) \longrightarrow \Diff(M)/\Diff_0(M).\]
\end{theorem}
Ontaneda \cite{Ontaneda} proved the above statement for the case of the double of a cusped hyperbolic manifold. His proof uses the fact that such a manifold admits a nonpositively curved metric and facts from CAT$(0)$ geometry. The proof we give does not use CAT$(0)$ geometry and thus covers the case of cusp-decomposable manifolds that do not admit a nonpositively curved metric, such as the case where the pieces in the cusp decomposition of $M$ is complex (or quarternionic or octonionic) hyperbolic.

\begin{remark}The requirement that gluing maps are affine diffeomorphisms is necessary for Theorem~\ref{Mostow} to be true. If the gluing maps are allowed to be any diffeomorphism, then the theorem is false. In \cite{AF}, Aravinda and Farrell prove that there exist gluing maps for the double of a hyperbolic cusped manifold that is homotopic to the identity map but the resulting space is not diffeomorphic to that obtained by the identity map.
\end{remark} 

The equation (\ref{*}) does not always split. Morita proved that if $M$ is a surface of genus greater than $4$, the group $\Out(\pi_1(M))$ does not lift to $\Diff(M)$ (see \cite[Theorem 4.21]{Morita}). However, for a cusp-decomposable manifold $M$ of dimension greater than $2$, the group $\Out(\pi_1(M))$ can be realized by a group of diffeomorphisms of $M$. Let 
\[\pi \colon \Diff(M) \longrightarrow \Diff(M)/\Diff_0(M)\] be the natural projection. 

\begin{theorem}\label{Out lifts}
Let $M$ be a cusp-decomposable manifold of dimension greater than $2$. There exists a homomorphism 
\[ \sigma \colon \Out(\pi_1(M)) \longrightarrow \Diff(M)\]
such that $\pi\circ\sigma = \rho$. Hence, if $\Out(\pi_1(M))$ is finite, then it can be realized as a group of isometries of $M$ with respect to some Riemannian metric on $M$.
\end{theorem}

By Mostow rigidity, if $M$ is a finite volume, complete, irreducible, locally symmetric, negatively curved manifold of dimension greater than $2$, then $\Out(\pi_1(M))$ is finite. In contrast, cusp-decomposable manifolds can exhibit a kind of non-rigidity. It follows from the proof of Theorem~\ref{Mostow} that if $M$ is a cusp-decomposable manifold of dimension greater than $2$, the group $\Out(\pi_1(M))$ is an extension of a finitely generated, torsion free, abelian group $\mathcal{T}(M)$ by a finite group $\mathcal{A}(M)$. (See Section~\ref{Out} for the definition of $\mathcal{T}(M)$ and $\mathcal{A}(M)$). The (infinite order) elements in $\mathcal{T}(M)$ can be realized by diffeomorphisms that are analogous to Dehn twists in surface topology.
\begin{theorem}\label{Out(pi_1(M))}
$\Out(\pi_1(M))$ is an extension of $\mathcal{T}(M)$ by $\mathcal{A}(M)$, i.e. the following sequence is exact.
\[1 \longrightarrow \mathcal{T}(M) \longrightarrow \Out(\pi_1(M)) \longrightarrow \mathcal{A}(M) \longrightarrow 1.\] 
The group $\Out(\pi_1(M))$ is infinite if and only if the fundamental group of one of the cusps in the cusp decomposition of $M$ has nontrivial center.
\end{theorem}

In the case where $\Out(\pi_1(M))$ is infinite, $\Out(\pi_1(M))$ cannot be realized by isometries of $M$ for any metric on $M$ (see \cite{Tam}), in contrast with the case of the Mostow Rigidity. As we will see, this happens to a lot of cusp-decomposable manifolds, such as the double of a one-cusp hyperbolic manifold whose cusp is diffeomorphic to the product of a torus and $\R$.
\newline
\newline
\textbf{Acknowledgements.} I would like to thank my advisor, Benson Farb, for suggesting this problem and for his guidance. I would like to thank Thomas Church and Matthew Day for helpful conversations. I would also like to thank Thomas Church, Spencer Dowdall and Benson Farb for commenting on earlier versions of this paper. 

\section{Cusp-decomposable manifolds: background and definition}
\label{sec:definitions}

Let $Y$ be a connected, noncompact, finite-volume, complete, locally symmetric,  negatively curved manifold of dimension $n \geq 3$. By the classification of  simply connected, complete, symmetric spaces of $\R$-rank $1$, the universal cover $\widetilde{Y}$ of $Y$ is isometric $\mathbb{KH}^n$, where $\mathbb{K}$ can be the real numbers $\R$, the complex numbers $\C$, the quarternions $\Q$ or the octonions $\A$ (in which case $n = 2$). By the thick-thin decomposition and a finiteness theorem for finite-volume, complete manifolds with pinched sectional curvature (see \cite{BGS}), $Y$ has finitely many ends or \emph{cusps}. Each cusp is diffeomorphic to $[0, \infty) \times S$ for some compact $(n-1)$--dim manifold $S$. The fundamental group of each cusp, or each \emph{cusp subgroup}, corresponds to a maximal subgroup of $\pi_1(Y)$ of parabolic isometries fixing a point on the boundary at infinity of $\widetilde{Y}$. The parametrization $[0,\infty) \times S$ of a cusp can be taken so that each cross section $a\times S$ is the quotient of a horosphere in $\widetilde{Y}$ by the corresponding cusp subgroup. 

By the Margulis lemma, cusp subgroups are virtually nilpotent. This also follows from the fact that the group of parabolic isometries of $\widetilde{Y}$ fixing a point on the boundary at infinity has a subgroup that is a $2$--step nilpotent Lie group that acts freely and transitively on the horosphere (\cite{GP}). In fact, $S$ with the induced Riemannian metric from $Y$ is a compact \emph{infra-nilmanifold}, that is, a compact quotient of nilpotent Lie group $H$ with left invariant metric by a torsion-free (necessarily cocompact) lattice of the group of isometries of $H$. 

If we delete the $(b, \infty) \times S$ part of each cusp $[0, \infty) \times S$ of $Y$, the resulting space is a compact manifold with boundary. We choose $b$ large enough so that the boundary components $b \times S$ of different cusps do not intersect. We require the boundary components of $Y$ lift to a horosphere in $\widetilde{Y}$. We call manifolds obtained this way \emph{bounded-cusp manifolds with horoboundary}. 

We say that a manifold $M$ is \emph{cusp-decomposable} if it is obtained by taking finitely many bounded-cusp manifolds  $Y_i$ with horoboundary and glue them along pairs of boundary components via affine diffeomorphisms. Each pair of boundary components that are glued together can belong to the same $Y_i$. The set of the spaces $Y_i$ with a gluing is called a (the) \textit{cusp-decomposition} of the manifolds. Each $Y_i$ is a \emph{piece} in the cusp decomposition of $M$. 

We give some more examples of cusp-decomposable manifolds apart from those in the Introduction. If $Y$ is a bounded-cusp manifold with two diffeomorphic horoboundary components  $b\times S$ and $b'\times S'$, then we can glue the two boundary components together by an affine diffeomorphism, which exists by the version of Bieberbach theorem for infra-nilmanifolds (see \cite{Auslander}). The resulting manifold is cusp-decomposable, and so is any finite cover of it.

The cusp decomposition of $M$ is analogous to the JSJ decomposition of 3-manifolds. That is, there exists a finite collection of embedded, codimension $1$, incompressible, infra-nil submanifolds  of $M$ such that if we cut $M$ along these submanifolds, the resulting space is a disjoint union of connected manifolds, each of which is diffeomorphic to a finite-volume, noncompact, complete, locally symmetric manifold of negative curvature. There is also a $2$-dimensional analog of the cusp decomposition with the pair of pants decomposition of a surface.

\section{Fundamental groups of cusp-decomposable manifolds}

\begin{figure}
\begin{center}
\includegraphics[height=100mm]{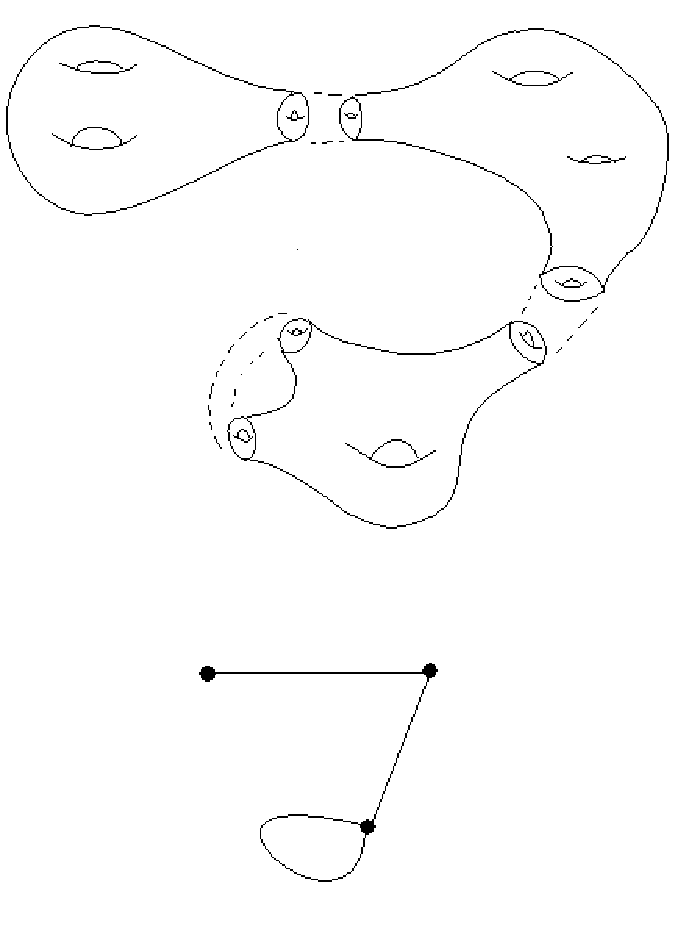}
\caption{A cusp-decomposable manifold and the corresponding underlying graph of the graph of groups structure of the fundamental group.}
\end{center}
\end{figure}

Let $M$ be a cusp-decomposable manifold of dimension $n >2$ and let $\{Y_i\}_{i \in I}$ be the pieces in the cusp decomposition of $M$. The cusp decomposition of $M$ gives $\pi_1(M)$ the structure of the fundamental group of a graph of groups $\mathcal{G}_M$ (see \cite{Serre} and \cite{Scott} for Bass-Serre theory). The vertex groups of $\mathcal{G}_M$ are the fundamental groups of the pieces $Y_i$ in the cusp-decomposition of $M$. The edge groups of $\mathcal{G}_M$ are the fundamental groups of the cusps of each piece, i.e. boundary components of each $Y_i$. The injective homomorphism from an edge group to a vertex group is the inclusion of the corresponding boundary component into $Y_i$. As we saw in the introduction, the edge groups are virtually nilpotent.

The following lemmas are facts about the fundamental group of a cusp-decomposable manifolds that we will use to prove Theorem~\ref{Mostow}.

\begin{lemma}\label{nilpotent subgroups}
Let $M$ be a cusp-decomposable manifold of dimension $n >2$. If $H$ is a virtually nilpotent subgroup of $\pi_1(M)$ that is not virtually $\Z$, then $H$ is contained in a conjugate of a vertex group of the graph of groups defined by the cusp decomposition of $M$.
\end{lemma}

Before proving Lemma~\ref{nilpotent subgroups}, we state a result in Bass-Serre theory (\cite{Serre}) which we will use in the proof.

\begin{theorem}[Serre]\label{nilpotent action}
Let $G$ be a finitely generated nilpotent group acting on a tree $T$ without inversion. Then either $G$ has a fixed point or there is an axis $L$ that is preserved by $G$, on which $G$ acts by translation by means of a nontrivial homomorphism $G \longrightarrow \Z$.
\end{theorem}

We will also describe a few properties of covering spaces of cusp-decomposable manifolds and say how their universal covers relate to their Bass-Serre tree. Let $\widehat{M}$ be a covering space of $M$ with covering map $p \colon \widehat{M} \longrightarrow M$. Since $M$ is obtained from gluing $Y_i$'s along their boundaries, each component of $p^{-1}(Y_i)$ in $\widehat{M}$ is a covering space of $Y_i$ for each $i$. Similarly, a component of $p^{-1}(\partial Y_i)$ in $\widehat{M}$ is a covering space of $\partial Y_i$ for each $i$. Thus, $\widehat{M}$ is the union of connected covering spaces of $Y_i$ glued along their boundaries.

By the above description of covering spaces of $M$, the universal cover $\widetilde{M}$ of $M$ is the union of universal covers of \emph{neutered spaces} $Y_i$ (which are complete simply connected locally symmetric negatively curved manifolds with a collection of disjoint horoballs removed) glued to each other along their horosphere boundaries. Each of the neutered spaces is a connected component of $p^{-1} (Y_i)$ for some $Y_i$ in the cusp decomposition on $M$. The underlying graph in the graph of groups defined by this decomposition of $\widetilde{M}$ is the Bass-Serre tree of the graph of groups $\mathcal{G}_M$. Two horospheres that belong to the same neutered space are called \emph{adjacent}.

\begin{proof}[Proof of Lemma~\ref{nilpotent subgroups}]
Let $K$ be a finite index nilpotent subgroup of $H$. By Theorem~\ref{nilpotent action}, the group $K$ acts on the Bass-Serre tree $T$ of $\pi_1(M)$ (with the structure of the graph of groups $\mathcal{G}_M$) either with a fixed point or preserving an axis $L$. We claim that the latter case cannot occur. Suppose $K$ preserves an axis $L$. Then the action of $K$ factors through a homomorphism $\varphi \colon K \longrightarrow \Z$, so $\Ker\varphi$ fixes the axis $L$. Now, we claim that any element $g \in \pi_1(M)$ that preserves a path in $T$ has to be the identity. Since $g$ fixes an edge, it belongs to a conjugate of an edge group. Also, $g$ corresponds to a parabolic isometry of a neutered space in the universal cover $\widetilde{M}$. Hence, if $g$ is nontrivial, it will act on the Bass-Serre tree of $\pi_1(M)$ by moving all the edges at one vertex (that is fixed) except one edge. So $g$ must be the identity element. Thus, $\Ker\varphi = 1$ and $K \cong \Z$, which contradicts the fact that $H$ is not virtually $\Z$. Therefore, $K$ fixes a vertex $a$ in $T$, since the action is without inversion. Thus, $K$ is contained in the stabilizer of $a$, which is a conjugate of a vertex group $G_{v}$ for some vertex $v$ in the Bass-Serre tree of $\pi_1(M)$.

Now we show that $H$ is contained in the same conjugate of $G_{v}$ as $K$. For each $h \in H$ that is not the identity element, there exists a positive integer $m$ such that $h^m \in K$, and thus $h^m$ fixes $v$. Since $h$ and $h^m$ commute, $h$ preserves the fixed point set of $h^m$, which, given the way $\pi_1(M)$ acts on $T$, is either $\{a\}$ or an edge containing $\{a\}$. In both cases, this means that $h$ fixes $a$ since $H$ acts on $T$ without inversions. Hence, $H$ is contained in the stabilizer of $a$ as well.   
\end{proof}

\begin{lemma}\label{Z^2 in cusp subgroup}
If $C$ is a cusp subgroup of a noncompact finite-volume locally symmetric negatively curved manifold $Y$ of dimension greater than $2$, then $C$ has a subgroup that is isomorphic to $\Z^2$.
\end{lemma}

\begin{proof}
Let $L$ be a finite-index nilpotent subgroup of $C$. Then the center $Z(L)$ of $L$ is nontrivial and thus infinite since $L$ is torsion-free. If $Z(L)$ contains a subgroup isomorphic to $\Z^2$, then we are done. If not, then $Z(L) \cong \Z$. Observe that $L$ is not virtually $\Z$ since $L$ is the fundamental group of a closed aspherical manifold of dimension at least $2$. Hence, $L/Z(L)$ is infinite. Let $a \in Z(L)$ and let $b \in L$ such that $b$ represents an infinite order element of $L/Z(L)$. Then the subgroup of $C$ generated by $a$ and $b$ is isomorphic to $\Z^2$.

\end{proof}

In order to prove Theorem~\ref{Mostow} we need to analyze isomorphisms between the fundamental groups of two cusp-decomposable manifolds. It turns out that an isomorphism preserves the graph of groups structure of these fundamental groups, as we will see in the next two lemmas.

\begin{lemma}\label{cusp subgroup}
Let $M$ and $N$ be two cusp-decomposable manifolds of dimension $n > 2$ and let $\phi \colon \pi_1(M) \longrightarrow \pi_1(N)$ be an isomorphism. If $G_e$ is an edge group in the graph of groups decomposition of $\pi_1(M)$, i.e. a cusp subgroup of $\pi_1(M)$,  then $\phi(G_e)$ is conjugate to an edge group of $\pi_1(N)$.
\end{lemma}

\begin{proof}
By Lemma~\ref{nilpotent subgroups}, the group $\phi(G_e)$ is a subgroup of a conjugate of a vertex group $G'_{v'}$ of $\pi_1(N)$. Let $L$ be a finite-index nilpotent subgroup of $\phi(G_e)$. By Lemma~\ref{Z^2 in cusp subgroup}, the group $\phi(G_e)$ has a subgroup $K$ that is isomorphic to $\Z^2$ . This implies that $K$ consists of only parabolic isometries and thus, is a subgroup of a cusp/edge subgroup $G'_{e'}$ of $G'_{v'}$.  So the fixed set of each nontrivial element of $K$ acting on the Bass Serre tree of $\pi_1(N)$ is an edge. 

Also, the proof of Lemma~\ref{Z^2 in cusp subgroup} shows that $K$ can be picked to contain a nontrivial element $c$ in the center of $L$. Since the fixed set of $c$ is an edge, each nontrivial element in $L$ fixes the same edge. So $L$ is contained in $G'_{e'}$. Now, for every nontrivial element $g \in \phi(G_e)$, there exists a positive integer $m$ such that $g^m \in L$ because $L$ is a finite index in $\phi(G_e)$. Since $g$ and $g^m$ commute, the fixed set of $g$ is the same as that of $g^m$. Hence, $g \in G'_{e'}$, so $\phi(G_e)$ is contained in  $G'_{e'}$. 
\end{proof}

\begin{lemma}\label{vertex subgroup}
Let $M$ and $N$ be two cusp-decomposable manifolds of dimension $n > 2$ and let $\phi \colon \pi_1(M) \longrightarrow \pi_1(N)$ be an isomorphism. If $G_v$ is a vertex group in the graph of groups decomposition of $\pi_1(M)$,  then $\phi(G_v)$ is conjugate to a vertex group of $\pi_1(N)$.
\end{lemma}

\begin{proof}[Proof of Lemma~\ref{vertex subgroup}]
By Lemma~\ref{cusp subgroup}, $\phi$ maps the cusp subgroups of $\pi_1(M)$ to conjugates of cusp subgroups of $\pi_1(N)$, and each cusp subgroup is the image of a conjugate of a cusp subgroup. Since $M$ and $N$ are aspherical with isomorphic fundamental groups, there exists a homotopy equivalence $f \colon M \longrightarrow N$ such that $f_* = \phi $. The lift $\widetilde{f} \colon \widetilde{M} \longrightarrow \widetilde{N}$ is a quasi-isometry since $M$ and $N$ are compact and $f$ is continuous.

Fix basepoints $x_0 \in M$ and $y_0 \in N$, and pick lifts $\widetilde{x}_0 \in \widetilde{M}$ and $\widetilde{y}_0 \in\widetilde{M}$, respectively. Without loss of generality, we can assume that $\widetilde{f}(\widetilde{x}_0) = \widetilde{y}_0$. Now, each horosphere in $M$ (or $N$) is the universal of a horoboundary component (which is a compact $(n-1)$-dim manifold) of a piece in the cusp decomposition of $M$ (or $N$). Since $\phi$ maps cusp subgroups to conjugates of cusp subgroups and since $\widetilde{f}$ is a $\phi$--equivariant quasi-isometry, the image of a horosphere under  $\widetilde{f}$ is a bounded distance from a horosphere (or a \emph{quasi-horosphere}), i.e. it is contained in a $d$--neighborhood of a horosphere. Hence, $\widetilde{f}$ map the horosphere containing $\widetilde{x}_0$ to the quasi-horosphere containing $\widetilde{y}_0$. Also, since $\phi$ has an inverse, every horosphere in $N$ is the image of a quasi-horosphere in $M$.

We want to prove that any two adjacent horospberes are mapped to two adjacent horospheres. Let $H_a$, $H_b$ be adjacent horospheres. Suppose that $\widetilde{f}(H_a)$ and $\widetilde{f}(H_b)$ are not adjacent quasi-horospheres, i.e. any path starting at a point in $\widetilde{f}(H_a)$ and ending at some point in $\widetilde{f}(H_b)$ crosses a quasi-horosphere $f(H_c)$ for some horosphere $H_c$. 

Since the $\widetilde{M}_i$ is negatively curved, the distance function to a horoball is strictly convex. Hence, for any $K >0$, we can pick a point $p \in H_a$ and a point $q \in H_b$ that have distance at least $K$ from $H_c$. Let $\gamma$ be a path connecting $p$ and $q$ such that the distance between $\gamma$ is at least $K$. Since $f$ is a quasi-isometry, the distance between $f(\gamma)$ and $f(H_c)$ is positive if $K$ is large enough, i.e. $f(\gamma)$ is a path connecting $f(H_a)$ and $f(H_b)$ that does not cross $f(H_c)$. But this is a contradiction to the assumption that $f(H_a)$ and $f(H_b)$ are separated by $f(H_c)$. So $f$ does not map adjacent horospheres to non-adjacent horospheres.

Since $f$ has an quasi-isometric inverse, $f$ maps non-adjacent quasi-horospheres to non-adjacent quasi-horospheres. This implies that $f_*$ maps each vertex group to a conjugate of a vertex group, which is what we want to prove.      
\end{proof}

\section{Rigidity of cusp-decomposable manifolds}

\begin{proof}[Proof of Theorem~\ref{Mostow}]
Let $\phi \colon \pi_1(M) \longrightarrow \pi_1(N)$ be an isomorphism.  We are going to construct a diffeomorphism from $M$ to $N$ by firstly defining diffeomorphisms on each piece in the cusp decomposition of $M$ and then gluing them together in such a way that the resulting diffeomorphism induces the isomorphism $\phi$. 

Let $M_i$ (and $N_j$ respectively), for $i$ in some index set $I$ (for $j$ in some index set $J$ respectively), be the pieces in the cusp decomposition of $M$ (and $N$ respectively). By Theorem~\ref{vertex subgroup}, $\phi$ defines a bijection $\alpha \colon I \longrightarrow J$ such that $\phi(\pi_1(M_i)) = \pi_1(N_{\alpha(i)})$ up to conjugation. By the Mostow Rigidity Theorem for finite-volume, complete, locally symmetric, negatively curved manifolds, the restriction of $\phi$ on each vertex group of $\pi_1(M)$ up to conjugation is induced by an isometry $f_i$ from $M_i$ to $N_{\alpha(i)}$.

At this point one may want to glue the isometries $f_i$ together and claim that it is the desired diffeomorphism from $M$ to $N$. However, there are two problems. One is that the gluing of $f_i$ at each pair of boundaries of $M_i$'s that are glued together might not be compatible with the gluing of $N_j$. The other problem is that if $M_i$ and $M_j$ are adjacent pieces, it may happen that ${f_i}_*$ and ${f_j}_*$ defines different isomorphisms (by a conjugate) on the fundamental group of the boundary shared by $M_i$ and $M_j$. (See Lemma~\ref{twist} below for such an example of $\phi$).

Hence, we need to do some modifications to $f_i$'s near the boundary of $M_i$ before gluing them together. We observe that the restriction of $f_i$ and $f_j$ to each pair of boundary components of $M_i$ and $M_j$ that are identified, say $S$, in the cusp decomposition of $M$ induces the same map from $\pi_1(S)$ to $\phi(\pi_1(S))$ up to conjugation. Since $S$ is aspherical, it follows that $f_i$ and $f_j$ are homotopic. Then by Theorem~\ref{infranil isotopy} below, we can modify the maps $f_i$ and $f_j$ by an isotopy of $S$ on a tubular neighborhood of the corresponding boundary component of $M_i$ and $M_j$ that is compatible the gluing of $N$ and the action of $\phi$ on fundamental groups. Let $f \colon M \longrightarrow N$ be the diffeomorphisms obtained by gluing the modified $f_i$. Then $f_* = \phi$.
\end{proof}

\begin{theorem}\label{infranil isotopy}
Let $S$ be a compact infra-nilmanifold whose universal cover is a $2$-step nilpotent Lie group $H$ with a left invariant metric. If $f \colon H \longrightarrow H$ is an $\pi_1(S)$-equivariant isometry that is $\pi_1(S)$-equivariant homotopic to the identity map of $H$, then $f$ is $\pi_1(S)$-equivariant isotopic to the identity map. 
\end{theorem}
Each boundary component of the pieces in the cusp decomposition of a cusp-decomposable manifold is a compact infra-nilmanifold whose universal cover is a horosphere in $\mathbb{KH}^n$, which is isometric to a ($2$-step) nilpotent Lie group $H$ with a left invariant metric. Therefore, the above theorem applies to the manifold $S$ in the proof of Theorem~\ref{Mostow}.
\begin{proof}
We write $S = H/\Gamma$ for some torsion free lattice $\Gamma$ in $\Isom(H)$. There is an explicit description of isometries of $H$ by the following theorem \cite{GW}. For a Lie group $G$ with Lie algebra $\mathfrak{g}$, we denote $\Aut(G)$ (respectively, $\Aut(\mathfrak{g})$) the group of automorphisms of $G$ (respectively, $\mathfrak{g}$). Then $\Aut(\mathfrak{g})$ is naturally isomorphic to $\Aut(G)$ if $G$ is simply connected. Given a left invariant Riemannian metric $g$ on $G$, we denote $\Stab(e)$ the stabilizer of the identity element $e \in G$ in $\Isom(G,g)$. We denote $\OO(g)$ the group of orthogonal linear transformations of $\mathfrak{g}$ with respect to the metric $g$ on $\mathfrak{g}$. 
\begin{theorem}[Gordon, Wilson]\label{Gordon, Wilson}
If $G$ is a simply connected upper triangular unimodular Riemannian Lie group, then 
\[\Stab(e) = \Aut(\mathfrak{g}) \cap \OO(g) ,\] 
Moreover, $G$ is normal in $\Isom(G)$, that is 
\[ \Isom(G) = G \rtimes (\Aut(\mathfrak{g}) \cap \OO(g)). \] 
\end{theorem}
Since $H$ is a ($2$--step) nilpotent Lie group, $H$ satisfies the assumption in Theorem~\ref{Gordon, Wilson}. We denote the Lie algebra of $H$ by $\mathfrak{h}$. We are going to use logarithmic coordinates on $H$ and express isometries of $H$ in terms of these coordinates. Since the exponential map is a global diffeomorphism from $\mathfrak{h}$ to $H$,  each element in $H$ is $\exp(u)$ for a unique $u \in \mathfrak{h}$. Using logarithmic coordinates means that we will write $u$ for $\exp(u)$. Let ``$*$'' denotes the group operation of $H$. We denote by ``$+$" the vector space addition operation on $\mathfrak{h}$. By the Baker-Campbell-Hausdorff formula, for $u, v \in H$, 
\[ u * v = u + v + \dfrac{1}{2}[u,v],\]
which means
\[\exp(u)\exp(v) = \exp(u + v + \dfrac{1}{2}[u,v])\]
in logarithmic coordinates.

We pick an orthonormal basis $\{x_1, x_2,\ldots, x_k, y_1, y_2,\ldots, y_l\}$ of the Lie algebra $\mathfrak{h}$ in which $\{ y_1, y_2, \ldots, y_l\}$ is an orthonomal basis of $Y : = [\h, \h]$. Obviously, any automorphism of $\h$ fixes $Y$. Let $X$ be the subspace of $H$ that is spanned by $\{ x_1, x_2, \ldots, x_k\}$. 

By Theorem~\ref{Gordon, Wilson}, an isometry $f$ of $H$ has the form $f(x) = \varphi(x) * b$, where $\varphi \in \Aut(\mathfrak{h}) \cap \OO(g) $, $b \in H$. Such an element $\varphi$ in $\Aut(\mathfrak{g}) \cap \OO(g)$ preserves the subspaces $X$ and $Y$, that is $\varphi$ has the form
\[ \varphi =  \left(\begin{matrix} A&0\\0&B  \end{matrix}\right),\]
where $A$ is some $k \times k$ orthogonal matrix and $B$ is some $l \times l$ orthogonal matrix. Note that not every matrix of the above form gives an isometry of $H$.

Since $S$ is an infra-nil manifold, by the version of Bieberbach's theorems for infra-nil manifolds (see \cite{Auslander}), $\Gamma$ has a finite index subgroup $\Gamma'$ that is lattice in $H$. Since $f$ is $\Gamma$- equivariantly homotopic to the identity, we have
\[ f \circ \gamma = \gamma \circ f\]
for all  $\gamma \in \Gamma$. Write $f$ as $f(v) =  \left(\begin{matrix} A&0\\0&B  \end{matrix}\right)v * b$ for $v \in H$. Let $\pi_X$ and $\pi_Y$ are projections onto $X$ and $Y$ respectively. Note that $\pi_X([u,v]) = 0$ for all $u, v \in H$. Thus, for $\gamma \in \Gamma'$ and $v \in \mathfrak{h}$, we have
\begin{align*} 
f \circ \gamma (v) &= A (\pi_X(v +\gamma + \dfrac{1}{2}[v,\gamma])) +  B (\pi_Y(v +\gamma + \dfrac{1}{2}[v,\gamma])) + b \\
 & + \dfrac{1}{2}[A (\pi_X(v +\gamma + \dfrac{1}{2}[v,\gamma])) +  B (\pi_Y(v +\gamma + \dfrac{1}{2}[v,\gamma])),b]\\ 
&= A (\pi_X(v +\gamma)) +  B (\pi_Y(v +\gamma + \dfrac{1}{2}[v,\gamma])) + b + \dfrac{1}{2}[A (\pi_X(v +\gamma)),b] 
\end{align*}
and
\begin{align*} \gamma \circ f (v) & = A(\pi_X(v ))+ B(\pi_Y(v)) + b + \dfrac{1}{2}[A(\pi_X(v ))+ B(\pi_Y(v)),b]  + \gamma\\
& + \dfrac{1}{2}[A(\pi_X(v )+ B(\pi_Y(v)) + b + \dfrac{1}{2}[A(\pi_X(v ))+ B(\pi_Y(v)),b], \gamma]\\
& = A(\pi_X(v ))+ B(\pi_Y(v)) + b + \dfrac{1}{2}[A(\pi_X(v )),b]  + \gamma + \dfrac{1}{2}[A(\pi_X(v )) + b, \gamma],
\end{align*}
The simplified form of the above two expressions is due to the fact that $H$ is $2$-step nilpotent, so the bracket with an element in $Y$ is $0$. 

Since $f \circ \gamma = \gamma \circ f$, 
and setting $ v = 0$, we have
\[ A(\pi_X(\gamma)) + B(\pi_Y(\gamma)) + \dfrac{1}{2}[A (\pi_X(\gamma)),b] = \gamma + \dfrac{1}{2}[b,\gamma].\]
Now apply $\pi_X$ to both sides of the equation, we get
\[ A(\pi_X(\gamma)) = \pi_X(\gamma).\]
Since $S$ is compact, by \cite[Theorem 5.1.6]{CG},  $\Gamma'$ contains a vector space basis of $\mathfrak{h}$. Thus, the above equation implies that $A$ is the identity matrix. Using this fact, we get

\[B(\pi_Y(\gamma)) + \dfrac{1}{2}[\pi_X(\gamma), b] = \pi_Y(\gamma) + \dfrac{1}{2}[b,\gamma]. \]
So
\[B(\pi_Y(\gamma))   = \pi_Y(\gamma) + [b,\gamma].\]
Since $\gamma$ takes values in a lattice of $Y$, it follows that $B$ is the identity matrix and $b$ is central in $H$. Hence, $f$ is a translation by $b$.  Now we take the straight line homotopy between $f$ and the identity map. It is not hard to show by writing the homotopy out in the logarithmic coordinates that it is $\Gamma$-equivariant. One simply check that for $ t \in [0,1]$,
\[(tf) \circ \gamma =  \gamma \circ (tf)\]
This follows easily given the fact that $f$ is a translation by $b$ (so $tf$ is a translation by $tb$) and the equation holds for $t = 1$.  Also, for each $t \in [0,1]$, the map $tf$ is injective. Therefore, $f$ is $\Gamma$-equivariantly isotopic to the identity map.
\end{proof}
The following corollary follows from the proof of Theorem~\ref{infranil isotopy}. We will need to use it in the proof of a theorem in Section~\ref{Out}.
\begin{corollary}\label{finite order isom}
Let $S$ be as in Theorem~\ref{infranil isotopy}. If $f \colon S \longrightarrow S$ is an isometry of order $p < \infty$ and $\widetilde{f}$ denotes the lift of $f$ to the universal cover $\widetilde{S}$, then $\widetilde{f}^p$ is the map defined by multiplication of an element in the center $\pi_1(S)$.
\end{corollary}

\section{$\Out(\pi_1(M))$: twists and turns}\label{Out}

If $M$ is finite-volume, complete, locally symmetric and negatively curved, then by the Mostow Rigidity Theorem, $\Out(\pi_1(M)) \cong \Isom(M, g_{loc})$, where $g_{loc}$ is the locally symmetric metric on $M$. This implies that $\Out(\pi_1(M))$ is finite since $\Isom(M)$ is finite. One might expect that if $M$ is cusp-decomposable, given the above theorem, $\Out(\pi_1(M))$ will also be finite. However, this is not true.  

For example, let $M$ be the double of a cusped locally symmetric negatively curved manifold. Let $M_1$ and $M_2$ be the two pieces in the cusp decomposition of $M$ and let $G_i = \pi_1(M_i)$. Then $\pi_1(M) = G_1*_CG_2$. Pick an element $c_0 \ne 1$ in the center of $C$. Let $\phi \colon \pi_1(M) \longrightarrow \pi_1(M)$ be induced by 
\[\phi(g) = 
\begin{cases}
&  g \; \qquad \qquad \text{if} \; g \in G_1 \\
& c_0gc_0^{-1} \qquad \text{if} \; g \in G_2.
\end{cases}\]
It is clear that $\phi$ extends to an automorphism of $G_1*_CG_2$ since $\phi$ is an automorphism when restricted to $G_1$ and $G_2$ and agrees on the intersection of $G_1$ and $G_2$. 

\begin{lemma}\label{twist}
Let $M$ and $\phi$ be as above. Then $\phi$ is an infinite order element of $\Out(\pi_1(M))$.
\end{lemma} 
\begin{proof}

For all $k \in \N$, $\phi^k$ is the identity on $G_1$ and conjugation by $c_0^k$ on $G_2$. Suppose $\phi^k$ is an inner automorphism of $G_1*_CG_2$ for some $k \in \N$. Then there exists $g \in G_1*_CG_2$ such that $\Conj(g)\circ\phi$ is the identity on $G_1*_CG_2$. This implies that for $g_1 \in G_1$, we have $g\phi^k(g_1)g^{-1} = g_1$. So $gg_1g^{-1} = g_1$ since $\phi^k$ is the identity on $G_1$. Thus, $g$ is in the centralizer of $G_1$ in $G_1*_CG_2$. 

We claim that $g$ is in the centralizer of $G_1$. Consider the action of $G_1*_CG_2$ on the Bass Serre tree $T$ of $G_1*_CG_2$. The fixed set of $G_1$ is a vertex $v$. Since $g$ commutes with every element of $G_1$, the fixed set of $G_1$ must be preserved by $g$. Hence, $g$ belongs to the stabilizer of $v$, which is $G_1$. Therefore, $g$ is in the centralizer of $G_1$ and thus $g =1$. This implies that $\phi$ is the identity homomorphism, which is a contradiction since $\phi$ acts nontrivially on $G_2$. Hence, $\phi$ is represents an infinite order element in $\Out(\pi_1(M))$.

\end{proof}

The automorphism $\phi$ above is induced by gluing the identity diffeomorphism on $M_1$ and a diffeomorphism on $M_2$ that ``twists" the boundary of $M_2$ around the loop $c_0$ as in the proof of Theorem~\ref{Mostow}. This is analogous to a Dehn twist in surface topology. For each loop $c$ in the center of $C$, we define a \emph{twist} around a loop $c$ to be a diffeomorphism constructed as above. By Lemma~\ref{twist}, twists induce infinite order elements of $\Out(\pi_1(M))$. 

Similarly, for each element in the center of a cusp subgroup of a cusp-decomposable manifold $M$, we define a corresponding \emph{twist} like above. The induced map of a twist on $\pi_1(M)$ has the same form as $\phi$, that is, up to conjugation, it is the identity on one vertex subgroup and conjugation by an element in the center of the edge subgroup on the other vertex subgroup. 

Let $\mathcal{T}(M)$ be the subgroup of $\Out(\pi_1(M))$ that is generated by twists. It is not hard to see that any two twists commute since either they are twists around loops in disjoint cusps or the loops they twist around commute since they are both in the center of the same cusp subgroup. Hence, $\mathcal{T}(M)$ is a torsion-free, finitely generated abelian group. Indeed, $\mathcal{T}(M)$ is isomorphic to the direct sum of the centers of the cusp subgroups of $M$. We observe that $\mathcal{T}(M)$ is a normal subgroup of $\Out(\pi_1(M))$. We claim that $\Out(\pi_1(M))/\mathcal{T}(M)$ is a finite group. We can actually write $\Out(\pi_1(M))$ as a group extension of $\mathcal{T}(M)$ by a group of isometries of a manifold. 

Let $\mathcal{M}$  be the disjoint union of all the complete locally symmetric spaces corresponding to the pieces $M_i$ in the cusp decomposition of $M$, i.e. they are the spaces before we delete their cusps.
By Theorem~\ref{vertex subgroup}, an element of $\Out(\pi_1(M))$ has to decend to an isomorphism between vertex groups up to conjugation, which, by the Mostow Rigidity Theorem, is induced by an isometry with respect to the complete, locally symmetric metric on each of the pieces. The decending map is a homomorphism from $\Out(\pi_1(M))$ to $\Isom(\mathcal{M})$. We call this induced map 
\[\eta \colon \Out(\pi_1(M))\longrightarrow \Isom(\mathcal{M}).\]
Let $\mathcal{A}(M) $ be the image of $\Out(\pi_1(M))$ under $\eta$. We call the elements of $\mathcal{A}(M)$ \emph{turns}. Then $\mathcal{A}(M)$ is the subgroup of $\Isom(\mathcal{M})$ of isometries of $\mathcal{M}$ whose restriction to the boundaries of each pair of cusps that are identified in the cusp decomposition of $M$ are homotopic with respect to the gluing. 

An isometry of $\mathcal{M}$ can permute the pieces $\mathcal{M}_i$'s in the cusp decomposition of $\mathcal{M}$. Thus, there is a homomorphism $\varphi$ from $\Isom(\mathcal{M})$ to the group of permutations of $k$ letters. Let $P$ be the image of $\varphi$. The kernel of $\varphi$ contains precisely those isometries of $\mathcal{M}$ that preserve each of the pieces in $\{\mathcal{M}_i\}$. Hence, $\Isom(\mathcal{M})$ has the structure of an extension of groups as follow. 
\[ 1 \longrightarrow \bigoplus_{i = 1}^k\Isom(\mathcal{M}_i) \longrightarrow \Isom(\mathcal{M}) \longrightarrow P \longrightarrow 1.\]
It follows that $\Isom(\mathcal{M})$ is finite since it has a finite index subgroup that is isomorphic to the direct sum of the isometry groups of the components of $\mathcal{M}$, which are finite. Thus, $\mathcal{A}(M)$ is finite. Now we prove Theorem~\ref{Out(pi_1(M))}.
\begin{proof}[Proof of Theorem \ref{Out(pi_1(M))}]
It suffices to prove that $\mathcal{T}(M) = \Ker\eta$. By definition, the action of a twist on the fundamental group of each of the pieces in the cusp decomposition of $M$ is the identity map up to conjugation by an element in the center of a cusp subgroup. By Mostow rigidity, the image of a twist under $\eta$ is the identity isometry.  Thus, $\mathcal{T}(M) \leq \Ker\eta$. 

Now if $\phi \in \Ker\eta$, then the isometry that $\phi$ induces on $\mathcal{M}$ is the identity map. Hence, for each $i$, the restriction of $\phi$ on $\pi_1(M_i)$ is the identity map in $\Out(\pi_1(M_i))$. Let $\phi_i$ be the restriction of $\phi$ to $\pi_1(M_i)$. If $M_i$ and $M_j$ are adjacent pieces that are glued together along $S$, then the restriction of $\phi_i$ to $\pi_1(S)$ is equal to that of $\phi_j$. This implies that $\phi_i$ and $\phi_j$ differ by a conjugation of an element in $\pi_1(S)$. Therefore, $\phi$ is a product of twists. So $\Ker\eta \leq \mathcal{T}(M)$.  
\end{proof}
Hence, $\Out(\pi_1(M))$ is an extension of a finitely generated abelian group by a finite group. Geometrically, Theorem~\ref{Out(pi_1(M))} says that an element $\Out(\pi_1(M))$ is a composition of twists and turns. A diffeomorphism that corresponds to an element $\phi$ in $\Out(\pi_1(M))$ can be obtained by firstly decomposing $M$ to its pieces $M_i$, then do the turns $\eta(\phi)$ on each $M_i$, and then glue these turns together with the right twists. Infinite order elements of $\Out(\pi_1(M))$ are those that has a power equal to a twist. The following theorem says that in some way finite order elements of $\Out(\pi_1(M))$ are more rigid.

\begin{theorem}\label{finite order elements}
Let $M$ and $\mathcal{M}$ be as above and $\phi \in \Out(\pi_1(M))$. If $\phi$ is nontrivial and has finite order and if the induced isometry on $\mathcal{M}$, i.e., $\eta(\phi)$, is the identity on some component of $\mathcal{M}$, then $\phi$ is the identity element.
\end{theorem} 
\begin{proof}
If $\eta(\phi)$ is the identity map on every component of $\mathcal{M}$, then by Theorem~\ref{Out(pi_1(M))}, $\phi$ is a twist and thus has infinite order. Hence, $\phi$ acts nontrivially on some component of $\mathcal{M}$. Hence, we can find two adjacent components $\mathcal{M}_1$ and $\mathcal{M}_2$ of $\mathcal{M}$ such that $\phi$ acts by the identity on $\mathcal{M}_1$ and nontrivially on $\mathcal{M}_2$. However, the restriction of the action of $\phi$ on a pair of boundary components of $\mathcal{M}_1$ and $\mathcal{M}_2$ that are glued together, say $S$, are homotopic (and thus isotopic). By Corollary~\ref{finite order isom}, a power of $\phi$ is a twist around a loop in $S$. But this implies that $\phi$ is infinite order in $\Out(\pi_1(M))$, which is a contradiction.
\end{proof}

\begin{proof}[Proof of Theorem \ref{Out lifts}]
We construct a homomorphism $\sigma \colon \Out(\pi_1(M)) \longrightarrow \Diff(M)$ such that $\pi\circ\sigma = \rho$. In this proof, instead of consider $M$, we consider a manifold $N$ that is diffeomorphic to $M$ that is obtained as follows. Take all the pieces $M_i$'s in the cusp decomposition of $M$. We can assume that the pieces $M_i$ are such that such that the number of isometry classes of the horo boundary components of $M_i$ is minimal. This is equivalent to saying that if $\varphi \in \Out(\pi_1(M))$, then for each piece $M_i$, the restriction of $\varphi$ to $\pi_1(M_i)$ to $\varphi(\pi_1(M_i))$ corresponds to an isometry of $M_i$ to $M_{\varphi(i)}$ (where $\varphi(i)$ denotes the index of the piece whose fundamental group is the image under $\varphi$ of a conjugate of $\pi_1(M_i)$). Now to each boundary component $\partial_jM_i$, we glue $S_{ij} = \partial_jM_i\times[0,1]$ along $\partial_jM_i \times 0$ by the identity map. Then $N$ is obtained by gluing corresponding $\partial_jM_i \times 1$'s together via the gluing maps given by the cusp decomposition of $M$. Clearly, $M$ and $N$ are diffeomorphic; and $N$ is obtained by gluing the pieces of $M$ together with tubes that are diffeomorphic to $\partial_jM_i \times [-1,1]$'s.

Now we define the homomorphism $\sigma$. Let $\varphi \in \Out(\pi_1(M))$. We define $\sigma(\varphi)$ to be the diffeomorphism $f$ of $N$ such that $f$ maps each $M_i$ isometrically to $M_{\varphi(i)}$. Any twists or turns will happen in the tube connecting the boundary components of the pieces as follows. The manifold $N$ is obtained by gluing the pieces of $M$ together with the tubes that are diffeomorphic to $\partial_jM_i \times [-1,1]$'s, each cross section of which has an affine structure inherited from that of the two boundary components $\partial_jM_i \times -1$ and $\partial_jM_i \times 1$ of the tube. This affine structure is induced by the affine structure of the boundary components of the pieces they get glued to. By the proof of Theorem \ref{infranil isotopy}, the restriction of $f$ to $\partial_jM_i \times (-1)$ and that to $\partial_jM_i \times 1$ differ by an element in the center of $H$, the universal cover of $\partial_jM_i$. We define $f$ to be the straight line homotopy (as in the proof of Theorem \ref{infranil isotopy}) on $\partial_jM_i \times [-1,1]$. 

We claim that $\sigma$ is a homomorphism. This is because for each $\varphi \in \Out(\pi_1(M))$, the map $\sigma(\varphi)$ is the unique diffeomorphism that induces $\varphi$ and that acts isometrically on each $M_i$ and the restriction of which is the straight line homotopy between the map on the two ends of each tube. The composition of such two maps is a map of the same type. By uniqueness, we see that $\pi\circ\sigma = \rho$. Hence, we have proved $\Out(\pi_1(M))$ lifts to $\Diff(M)$.

If $\Out(\pi_1(M))$ is finite, by the above, it is realized by a group of diffeomorphisms $F$ of $M$. Pick any Riemannian metric $g_0$ on $M$ and average this metric by $F$ to get a metric $g$. Then $F$ is a group is a group of isometries of $M$. 
\end{proof}

\bibliographystyle{amsplain}
\bibliography{bibliography}


\end{document}